\documentclass[a4paper,reqno]{amsart}
\usepackage{amssymb,amscd}
\usepackage{color}

   \renewcommand{\a}{\alpha}

  \newcommand{\w}{\omega}

   \newcommand{\R}{\mathbb{R}}
   
      \newcommand{\T}{\mathbb{T}}

   \newcommand{\cF}{\mathcal{F}}

\newcommand{\SL}{{\mathcal{L}}}

\newtheorem{theorem}{{\bf Theorem}}
\newtheorem{corollary}{{\bf Corollary}}
\newtheorem{definition}{ {\bf Definition}}

\newtheorem{example}{Example}
\newtheorem{remark}{Remark}

\begin{document}

\author{David Mart\'inez Torres} \address{ Departamento de Matem\'{a}tica
PUC Rio de Janeiro \\  Rua Marqu\^{e}s de S\~{a}o Vicente, 225 \\ G\'{a}vea, Rio de Janeiro - RJ, 22451-900
\\ Brazil}\email{dfmtorres@gmail.com}
\author{Eva Miranda}
\address{Laboratory of Geometry and Dynamical Systems \\ Department of Mathematics\\
 Universitat Polit\`{e}cnica de Catalunya and BGSMath, EPSEB\\
Avinguda del Doctor Mara\~{n}\'{o}n 44-50\\
08028, Barcelona, Spain  \newline  and \\ \newline IMCEE, Observatoire de Paris \\ 77, avenue Denfert Rochereau
75014 Paris, France}
\email{eva.miranda@upc.edu}
\thanks{Eva Miranda is supported by the Catalan Institution for Research and Advanced Studies via an ICREA Academia 2016 Prize, by a \emph{Chaire d'Excellence} de la Fondation Sciences Math\'{e}matiques de Paris, and partially supported by the  Ministerio de Econom\'{\i}a y Competitividad project with reference MTM2015-69135-P/FEDER and by the Generalitat de Catalunya project with reference 2014SGR634. This work is supported by a public grant overseen by the French National Research Agency (ANR) as part of the \emph{\lq\lq Investissements d'Avenir"} program (reference: ANR-10-LABX-0098).
 Part of this research has been supported by an Ajut mobilitat 2017-EPSEB. We thank the Centre de Recerca Matem\`{a}tica CRM-Barcelona for their hospitality in July 2017 when this project has been finalized.}

\title{ Zeroth Poisson homology, foliated cohomology and perfect Poisson manifolds }
 \maketitle

\begin{abstract} We prove that for regular Poisson manifolds, the zeroth homology group is isomorphic
to the top foliated cohomology group and we give some  applications. In particular,
we show that for regular unimodular Poisson manifolds top Poisson and foliated cohomology groups are isomorphic. Inspired by the symplectic setting, we define what is a perfect Poisson manifold. We  use these Poisson homology computations to provide families of perfect Poisson manifolds.
\end{abstract}
\section{Introduction}

The description of geometric and algebraic properties of infinite dimensional groups and algebras is an important mathematical problem.
Some of these groups are linked to physical problems,  giving an additional motivation to pursue their study.
An important example is that of the group of volume preserving transformations of a manifold. By Arnold's
work \cite{ar} their Riemannian geometry is deeply related to fluid mechanics. Another relevant example is the group of Hamiltonian diffeomorphisms of phase spaces. Since Calabi's proof on the perfectness
of the Lie algebra of Hamiltonian vector fields on a symplectic manifold \cite{calabi}, much has been
done on the study of the group of Hamiltonian transformations of a symplectic manifold \cite{banyaga,ono}.
However,  little
is known on the structure of the group of Hamiltonian transformations of a Poisson manifold.

In this paper we address Calabi's question for a regular Poisson manifold, that is, the study
of the subalgebra of commutators of its Poisson algebra. It is well-known that
such subalgebra is captured by the zeroth Poisson homology of the Poisson manifold.  Poisson homology and cohomology
groups are the main algebraic invariants of a Poisson manifold. However, it is in general very difficult to
discuss their structure, and even less to provide explicit computations
(\cite{vaisman}, \cite{lichnerowicz},\cite{etingof}, \cite{pichereau},
 \cite{guimipi}). The characteristic foliation
of a regular Poisson manifold is the (regular) foliation integrating the Hamiltonian directions. Thus,
properties of the the commutator subalgebra are necessarily tied to properties of the characteristic foliation.
Our main result in this paper follows closely Calabi's original arguments to prove that the zeroth cohomology group of a regular Poisson manifold is
entirely controlled by the behaviour of the characteristic foliation:
 \begin{theorem}\label{thm:1} Let $(M,\pi)$ be a Poisson regular Poisson manifold of rank $2n$,
 and let $\cF_\pi$ denote its characteristic foliation.
  Then there is a canonical isomorphism of vector spaces:
  \[H_0^\pi(M)\longrightarrow H^{2n}(\cF_\pi)\]
 \end{theorem}

Our result provides a de Rham approach of the zeroth Poisson homology group which
is more amenable to explicit computations.  Also under the unimodularity assumption\footnote{Poisson manifolds admitting an invariant measure with respect to Hamiltonian vector fields are called unimodular (unimodularity can be captured using Poisson cohomology).} it
identifies the top Poisson and foliated cohomology groups:
\begin{corollary}\label{cor:1} Let $(M,\pi)$ be an orientable unimodular regular Poisson manifold of dimension $m$ and rank $2n$. Then
 then there is an isomorphism of cohomology groups:
 \[H^m_\pi(M)\rightarrow H^{2n}(\cF_\pi)\]
\end{corollary}

It has been recently shown \cite{DKM,BDMOP} that unimodular Poisson manifolds appear as geometrical model for the critical set (collision) of several problems in celestial mechanics which were largely studied by Arnold \cite{ar2} (see also \cite{chenciner}).

The structure of the paper is as follows: in Section 1 we prove Theorem \ref{thm:1}.
In Section 2 we discuss applications of Theorem \ref{thm:1} and exhibit various examples where the group $H_0^\pi(M)$ is computed. This provides families of examples of perfect Poisson manifolds.

\section{Commutators and foliated cohomology}

\textbf{Notation:} In this paper $H^k_\pi(M)$ will denote the degree  $k$ the Poisson cohomology group
\cite{lichnerowicz}, $H_k^\pi(M)$ with denote the degree $k$-homology group as defined by Brylinski \cite{brylinksi},
and $H^{k}(\cF_\pi)$ will denote the degree $k$ foliated cohomology group of $(M,\cF_\pi)$.
When we refer to top foliated cohomology groups we will mean $H^{m}(\cF_\pi)$,  where $m$ the rank of the (characteristic)
foliation.

We begin with the proof of theorem \ref{thm:1}:

\begin{proof}
We regard the Poisson structure $\pi$ as a closed, non-degenerate foliated 2-form $\w_{\cF_\pi}\in \Omega^2(\cF_\pi)$,
where closedness is with respect to the foliated de Rham differential $d_{\cF_{\pi}}$.

Recall that by definition:
\[H^\pi_0(M):=\frac{C^\infty(M)}{\{C^\infty(M),C^\infty(M)\}}.\]

We define the map:
\begin{eqnarray}\nonumber
 \phi: C^\infty(M) &\longrightarrow & \Omega^{2n}(\cF_\pi)\\
       f &\longmapsto &f  \frac{\w^n_{\cF_\pi}}{n!}.
\end{eqnarray}
Hence to prove the theorem we must show the equality
\[\phi(\{C^\infty(M),C^\infty(M)\})= d_{\cF_{\pi}}(\Omega^{2n-1}(\cF_\pi)).\]

To check the inclusion $\phi(\{C^\infty(M),C^\infty(M)\})\subset d_{\cF_{\pi}}(\Omega^{2n-1}(\cF_\pi))$ we use Leibniz's rule
(which holds in the foliated setting whenever they involve vector fields tangent to the foliation):
\[L_{X_f}(g\frac{\w^n_{\cF_\pi}}{n!})=(L_{X_f}g)\frac{\w^n_{\cF_\pi}}{n!}+g L_{X_f}\frac{\w^n_{\cF_\pi}}{n!}.\]
Since the leafwise Liouville volume form  $\frac{\w^n_{\cF_\pi}}{n!}$ is Hamiltonian invariant the following equality holds
\[\{f,g\}\frac{\w^n_{\cF_\pi}}{n!}=L_{X_f}(g\frac{\w^n_{\cF_\pi}}{n!})=d_{\cF_\pi}(i_{X_f}g\frac{\w^n_{\cF_\pi}}{n!}),\]
as we wanted to prove.

To prove the other inclusion we consider $\gamma\in d_{\cF_{\pi}}(\Omega^{2n-1}(\cF_\pi))$. Let us choose a locally finite open cover $U_i$, $i\in I$, such that each $U_i$ is in the domain of Weinstein local coordinates (as the splitting theorem guarantees \cite{weinstein3}) and let us
fix a partition of unity $\beta_j$, $j\in J$, subordinated to this open cover.

By setting $\gamma_j:=\beta_j\gamma$, then it is clear
that our problem reduces to showing that $d_{\cF_\pi}\gamma_j$ can be written as a sum of commutators supported in $U_{i(j)}$.

If we let $x_1,\dots,x_{2n}\in C^\infty(U_{i(j)})$ be the (symplectic) Weinstein local coordinates, then
in $U_{i(j)}$ we can write:
\[{\gamma_j}=\sum_{i=1}^{2n} h_i dx_1\wedge\cdots \wedge d\hat{x}_i\wedge \cdots \wedge dx_{2n},\quad h_i\in C^\infty(U_{i(j)}).\]
Therefore:
\[d_{\cF_{\pi}}\gamma_j=\sum_{i=1}^{2n}(-1)^{i-1}n!\frac{\partial h_i}{\partial x_i}\frac{\w^n_{\cF_\pi}}{n!}\]
We define $h_i:=\beta x_i$, where $\beta$ is supported in $U_{i(j)}$ and equals 1 in the support of $\gamma_j$. Then
\[\sum_{i=1}^n (-1)^{2i-1}n!\{h_{2i_1},h_{2i}\}+\sum_{i=1}^n (-1)^{2i-1}n!\{h_{2i},h_{2i-1}\}=\sum_{i=1}^{2n}(-1)^{i-1}n!\frac{\partial h_i}{\partial x_i},\]
and this finishes the proof of the theorem.

\end{proof}

Theorem \ref{thm:1} has an interesting consequence: let $(M,\cF)$ be a foliated manifold and consider
the set of Poisson structures whose characteristic foliation is $\cF$.  Observe that this
 subset of the leafwise closed and non-degenerate 2-forms  is open in the compact open topology. The interesting observation is that
any Poisson structure on this set has zeroth Poisson homology group canonically isomorphic to $H^{\mathrm{top}}(\cF)$ so the zeroth Poisson homology group only depends on the characteristic foliation.

\begin{remark} Our strategy to prove Theorem \ref{thm:1} is analogous to that of Calabi \cite{calabi}.
Lichnerowicz \cite{lichnerowicz} also used the same approach to prove that if $U\subset M$
 is contractible and contained in the domain of Weinstein coordinates, then
 ${C^\infty(U)}={\{C^\infty(U),C^\infty(U)\}}$. However, his proof is entirely local and does not make
 any reference to the foliated de Rham differential.

 Our result can also be derived from a spectral sequence
 argument to compute Poisson homology  \cite{vaisman}. We believe our proof is more transparent and direct.

\end{remark}

\begin{remark} This isomorphism generalizes Brylinski's isomorphism \cite{brylinksi} at any degrees $H_m^{\pi}(M^{2n})\cong H^{2n-m}(M^{2n})$   for symplectic manifolds.
\end{remark}

\begin{remark} A different connection between foliated cohomology and Poisson cohomology already appears in \cite{weinstein}
where  the modular class (see subsection \ref{ssec:unimod}) is related to the Reeb class of the symplectic foliation.
 This is made more precise by \cite{boucetta} where an injection $H^1(\mathcal F)\to H^1_{\pi}(M)$ is exhibited.

 \end{remark}

\section{Applications}

In this section we describe several applications of the main theorem.

\subsection{Unimodular Poisson manifolds}\label{ssec:unimod}

For a given volume form $\Omega$ on an oriented Poisson manifold $M$ the associated
modular vector field is defined  \cite{weinstein} as the following derivation:
 \[
C^\infty(M) \to \R : \, f\mapsto \frac{\SL_{X_f}\Omega}{\Omega}.
\]
It is a Poisson vector field and it preserves the volume form.
Given two different choices of volume form $\Omega$, the resulting modular vector fields  differ by a Hamiltonian
vector field thus defining the same class in the first Poisson cohomology group. This class is known as the \textbf{modular class} of the Poisson manifold. A Poisson manifold is called  \textbf{unimodular} if its modular class vanishes.

We recall from \cite{weinstein2} the duality between Poisson cohomology and homology for this class of Poisson manifolds.

\begin{theorem}[Evens-Lu-Weinstein \cite{weinstein2}]\label{thm:2} Let $(M,\pi)$ be an $m$-dimensional
orientable unimodular regular Poisson manifold then,
\[
H^k_{\pi}(M)\cong H^{m-k}_{\pi}(M).
\]
\end{theorem}

As a consequence of  Theorems \ref{thm:1} and \ref{thm:2} we obtain the following theorem
which relates Poisson and foliated cohomology of top degree:

\begin{corollary}\label{cor:1} Let $(M,\pi)$ be an orientable unimodular regular Poisson manifold of dimension $m$
and rank $2n$. There is an isomorphism of cohomology groups:
\begin{equation}\label{eq:isouni}
 H^{m}_\pi(M)\rightarrow H^{2n}(\cF_\pi).
  \end{equation}

\end{corollary}

As we remarked in the introduction, for a general Poisson manifold
little can be said on its
of Poisson (co)homology groups. However, the are some important types of Poisson structures
for which much more can be said, and, in particular, isomorphism (\ref{eq:isouni}) holds true. As these
types of Poisson manifolds are unimodular, isomorphism  (\ref{eq:isouni}) can be obtained by simply applying
Corollary \ref{cor:1}.
We illustrate this with two examples:

\subsubsection{\textbf{Example 1: Cosymplectic structures}} Cosymplectic manifolds have widely been studied related to different problems in Differential Geometry and Topology (see for instance de survey \cite{survey} and \cite{Marisa}). They are also ubiquitous in Poisson Geometry since their geometrical data determine a codimension-one symplectic foliation. Examples of cosymplectic manifolds can be found in \cite{guimipi0} and \cite{guimipi} since the critical set of a $b$-symplectic manifold has naturally cosymplectic structures associated to them. In particular so are the collision set \cite{DKM} of the geometrical models for the n-body problem and other problems in celestial mechanics \cite{BDMOP} as studied by Arnold \cite{ar2}.

A {\textbf{ cosymplectic structure}} on a manifold $M^{2n+1}$ is given by a pair of
closed forms $(\theta,\eta) \in \Omega^1(M) \times \Omega^2(M)$
for which $\theta\wedge \eta^{n}$ is a volume form on $M$.  Cosymplectic manifolds are naturally  endowed with  a
corank-one regular Poisson structure whose
characteristic foliation integrated the kernel of $\alpha$.

\begin{theorem}[Osorno, \cite{boris}]\label{thm:boris} The Poisson cohomology of the Poisson
structure $\pi$ associated to a cosymplectic manifold can be given in terms of the foliated cohomology as follows:
$$H^{k}_{\pi}(M)\simeq H^k(\mathcal{F_\pi})\oplus  H^{k-1}(\mathcal{F_\pi}).$$
\end{theorem}

For the top Poisson cohomology group the above theorem gives: $H^{2n+1}_\pi(M)\cong H^{2n}(\cF)$. Since
for a cosymplectic manifold $\a\wedge \eta^n$ is a Hamiltonian invariant volume form, we can obtain
this isomorphism by simply applying Corollary \ref{cor:1}.

\subsubsection{\textbf{Example 2: Poisson manifolds of s-proper type}}
A Poisson manifold is said to be of compact type the Lie algebroid $A=(T^*M,\sharp)$
associated to the Poisson structure $\pi$ with anchor map $\sharp: T^*M\to TM,\sharp(\alpha):= \pi(\alpha, \cdot)$
integrates to a compact-like \emph{symplectic Lie groupoid} \cite{crainicetal}. There exist several
natural compact type conditions one of which
defines the so-called Poisson manifolds of source proper type (they are integrated by some source proper
symplectic groupoid).

Hodge-type isomorphisms hold for Poisson manifolds of source proper type \cite{crainicetal}:
$$H^k_{\pi}(M)\cong H^{top-k}_{\Pi}(M).$$
As Poisson manifolds of compact types are unimodular \cite{crainicetal}, one deduces:
\begin{theorem}[Crainic-Fernandes-Martinez,\cite{crainicetal}] \label{thm:3} For  orientable Poisson manifolds
of source proper type there exist isomorphisms:
$$H_k^{\pi}(M)\cong H^k_{\pi}(M).$$
\end{theorem}

Theorem \ref{thm:3} in degree zero says
that $H_0^\pi(M)$ is isomorphic to the Casimirs of $(M,\pi)$. As the characteristic foliation
of a Poisson manifold of source proper type is compact and has finite holonomy  \cite{crainicetal}, the
Casimirs are isomorphic to the foliated cohomology group of top degree, this giving (\ref{eq:isouni}).
Once more, unimodularity suffices to obtain  Theorem \ref{thm:3} (for all compact types).

\subsection{Perfect Poisson manifolds}

We would like to finish this note by going back to the discussion in the introduction. Calabi's original motivation
was proving that the Lie algebra group of symplectomorphisms is perfect.

If we consider the Lie algebra of Hamiltonian vector fields we may want to see that it is perfect in a direct way.
 Because it is easier to study Poisson brackets of functions rather than Lie brackets of vector fieldsl we will compare the set of Hamiltonian vector fields with the set of functions in the symplectic setting.

For a symplectic manifold  $(M,\omega)$  the  Lie algebra of Hamiltonian vector fields  is perfect.
 The mapping that sends a function to its Hamiltonian vector field, $f\mapsto X_f$, defines
the short exact sequence of Lie algebras,

\[ 0\to \R\to C^\infty(M)\to \mathrm{ham}(M,\omega)\to 0.\]
It follows that $\mathrm{ham}(M,\omega)$ is perfect if and only if
\[C^\infty(M)=\{C^\infty(M),C^\infty(M)\}+\R.\]

We may use the quotient $\frac{C^\infty(M)}{\{C^\infty(M),C^\infty(M)\}}$ as a first approach to this equality. This quotient coincides with the zeroth Poisson homology $H_0^\pi(M)$.
From Brylinki \cite{brylinksi} the zeroth Poisson homology of a symplectic manifold is isomorphic to the top de Rham cohomology class which, in turn, is of dimension 1 in the compact case. Thus a symplectic manifold satisfies the equality above and a symplectic manifold is perfect.

 For a  Poisson manifold $(M,\pi)$
the analogous question is whether the Lie algebra of Hamiltonian vector fields associated to a given Poisson structure $\mathrm{ham}(M,\pi)$ is perfect.
From the discussion above in the symplectic setting we adopt the following definition of \textbf{perfect Poisson manifold.}
\begin{definition} A Poisson manifold is perfect if and only if, the following equality holds:
\[C^\infty(M)=\{C^\infty(M),C^\infty(M)\}+\R.\]

\end{definition}

Thus perfectness of the set of Hamiltonian vector fields is connected to the zeroth-Poisson homology group $H_0^\Pi(M)$ of a Poisson manifold.

In particular the discussion above proves.

\begin{theorem} If $H_0^\pi(M)=0$ then $(M, \pi)$ is perfect.

\end{theorem}

 This provides a wide class of examples. Let us construct one of them:

\subsubsection{Example: A corank one Poisson manifold which is perfect and not unimodular}\label{section:boucetta}

We consider the $3$-dimensional Poisson manifold with symplectic foliation defined as follows:

\begin{itemize}

\item Following  \cite{elkacimi} we denote by $\T^3_A$ the mapping torus associated to the diffeomorphism on $\T^2$ given by a matrix in $SL(2,\mathbb Z)$. In other words, on the product $\T^2\times \mathbb R$ consider the action: $(m,t)\to (A(m), t+1)$ for $A\in SL(2,\mathbb R)$ and we identify $(m,t)$ with $(A(m), t+1)$. For this construction we assume  $A$ to be a hyperbolic matrix i.e., $\vert Tr(A)\vert >2$.

   \item  In $\T^3_A$   consider the $2$-dimensional foliation as follows: For an hyperbolic element $A$ there are two real  eigenvalues with irrational slope. We denote by $\lambda$ one of the eigenvalues and by $\alpha$ the irrational slope. The foliation by lines of irrational slope $\alpha$ on the torus $\T^2$ is invariant by the action described above and thus descends to the quotient $\T^3_A$.
    Denote by $\mathcal F$ the foliation generated in this way on $\T^3_A$.
\item Observe that $\T^3_A$  is naturally endowed with a Poisson structure whose induced symplectic structure on each of the leaves is given by the area form.

    \item In \cite{elkacimi} the foliated cohomology of $\T^3_A$ is computed. In particular:

    \begin{theorem}[El Kacimi, \cite{elkacimi}]\label{thm:4}
$    H^2(\mathcal F)=0.$
    \end{theorem}
    \end{itemize}

As a consequence of  theorems \ref{thm:1} and \ref{thm:4} we obtain the following Poisson homology computation,

\begin{corollary} For the Poisson manifold $\T^3_A$ described above,
$H_0^\pi(\T^3_A)=0$ and, in particular, the Poisson manifold $\T^3_A$ is perfect.

\end{corollary}

\begin{remark} As observed already in \cite{elkacimi} since the top foliated cohomology group of this foliation does not vanish then it does not admit any transverse density. In particular,   the codimension-one foliation constructed above is not unimodular.

\end{remark}

For unimodular Poisson manifolds the restrictions on the zeroth homology group for the Poisson manifold to be perfect are weaker.

Recall that an unimodular Poisson manifold admits a Hamiltonian invariant volume. In \cite{weinstein} it is shown that that for $\mu$ any  Hamiltonian invariant volume
form and $C^{\infty}_\mu(X)$ the hyperplane of $C^{\infty}(X)$ of zero mean
functions with respect to $\mu$, one has
\begin{equation}\label{eq:zeromean}
\{C^{\infty}(X),C^{\infty}(X)\}\subset C^{\infty}_\mu(X)
\end{equation}
In particular, for unimodular Poisson manifolds we obtain,

\begin{theorem}\label{thm:unimod} If a Poisson manifold $(M,\pi)$  is unimodular and satisfies $\dim H_0^\pi(M)=0,1$   then it is perfect.

\end{theorem}
From this theorem we may recover that any compact symplectic manifold is perfect (as we already proved above) because it is unimodular and $\dim H_0^\pi(M)=1$.
\begin{example}\label{ex:x}

In particular, as an application of theorem \ref{thm:unimod} a Poisson manifold $(M_1, \pi_1)$ given by any cosymplectic manifold in dimension 3 with a compact symplectic leaf which is a 2-sphere $S^2$  (and thus all of them are, \cite{guimipi}) is perfect.

\end{example}
We may produce product-type examples (and counter-examples) from the ones above using a K\"{u}nneth-type formula for foliated cohomology.

The following statement was proved for sheaf cohomology in the context of Geometric Quantization \cite{mirandapresas}. When applied to the limit case $\omega\to 0$ it yields the following  K\"{u}nneth formula for foliated cohomology of the product foliated manifolds $(M_1,\mathcal{F}_1)\times (M_1,\mathcal{F}_2)$.

\begin{theorem}[Miranda-Presas \cite{mirandapresas}] \label{thm:Kun}
     There is an isomorphism
$$ H^n(\mathcal{F}_1\times \mathcal{F}_2 ) \cong \bigoplus_{p+q=n} H^p(\mathcal{F}_1) \otimes H^q(\mathcal{F}_2), $$
\noindent whenever the foliated cohomology associated to $\mathcal{F}_1$ has finite dimension, $M_1$ is compact  and $M_2$ admits a \emph{good} covering.
\end{theorem}

From this theorem we deduce:

\begin{example} The product of two different cosymplectic manifolds $(M_1, \pi_1)$ in example \ref{ex:x} is perfect.(This includes the product of cosymplectic manifolds with symplectic manifolds).

\end{example}

\begin{example} The product of  $(\T^3_A, \pi)$ constructed in section \ref{section:boucetta} with any symplectic manifold is  perfect.
\end{example}

\begin{example} The product of  any cosymplectic manifold with the Poisson manifold $(\T^3_A, \pi)$ constructed in section \ref{section:boucetta} is  perfect.

\end{example}

\end{document}